\documentclass[11pt]{article}
\usepackage{amsmath}
\usepackage{amsfonts}
\newtheorem{thm}{Theorem}[section]
\newtheorem{lem}[thm]{Lemma}

\newtheorem{defn}[thm]{Definition}
\newenvironment{defn-new}{\begin{defn} \em}{\end{defn}}
\newtheorem{rem}[thm]{Remark}
\newenvironment{rem-new}{\begin{rem} \em}{\end{rem}}
\newtheorem{ex}[thm]{Example}
\newenvironment{ex-new}{\begin{ex} \em}{\end{ex}}

\newenvironment{notation-new}{\begin{rem} \em}{\end{rem}}

\newenvironment{agr-new}{\begin{rem} \em}{\end{rem}}

\makeatletter \@addtoreset{equation}{section} \makeatother

\makeatletter \@addtoreset{figure}{section} \makeatother

\begin{document}
\begin{center}
    {\Large\textbf {DDVV conjecture for Riemannian maps from quaternionic space forms}}\\[0pt]

 {\bf Kirti Gupta}\footnote{
Department of Mathematics \& Statistics, Dr. Harisingh Gour Vishwavidyalaya,
Sagar-470 003, M.P. INDIA \newline
Email: guptakirti905@gmail.com}, {\bf Punam Gupta}\footnote{
School of Mathematics, Devi Ahilya Vishwavidyalaya, Indore-452 001, M.P. 
INDIA\newline
Email: punam2101@gmail.com}  and {\bf R.K. Gangele}\footnote{
Department of Mathematics \& Statistics, Dr. Harisingh Gour Vishwavidyalaya,
Sagar-470 003, M.P. INDIA \newline
Email: rkgangele23@gmail.com} 
\end{center}

\noindent {\bf Abstract:}  
In this paper, we investigate the DDVV-type inequality for Riemannian maps from quaternionic space forms to Riemannian manifolds. We also discuss the equality case of the derived inequality with application.


\section{Introduction}
In the field of geometry, inequalities serve to link intrinsic and extrinsic invariants, playing a crucial role in defining the characteristics of geometric objects. Among the most significant results in this area is the DDVV conjecture.

The DDVV conjecture has its roots in the earlier work of Pierre Wintgen \cite{Wintgen}, who in 1979 proved a seminal inequality for surfaces immersed in a $4$-dimensional Euclidean space. This result was later generalized by Guadalupe and Rodriguez \cite{Guad} to surfaces in arbitrary space forms. The DDVV conjecture, formally proposed in 1999 by De Smet, Dillen, Verstraelen, and Vrancken \cite{DDVV}, extends this foundational work to higher dimensions and arbitrary codimension.

A DDVV-type inequality is an estimate of the form
\begin{equation}
\sum_{r,s=1}^m ||A_r, A_s||^2 \leq c \left(\sum_{r=1}^m ||A_r||^2\right)^2.
\end{equation}
This inequality is considered for a certain type of $n\times n$ matrices $A_1,\ldots,A_m,$ where $[A,B]=AB-BA$ denotes the commutator. The term $||A_r||^2
=tr(AA^{\ast})$ represents the squared Frobenius norm, where $A^{\ast}$
  is the conjugate transpose, and $c$ is a non-negative constant. The DDVV-type inequality originates from the normal scalar curvature conjecture (also known as the DDVV conjecture) in submanifold geometry.

In 1999, De-Dillen-Verstraelen-Vrancken \cite{DDVV} proposed the normal scalar curvature conjecture:

\noindent\textbf{Conjecture:} (DDVV conjecture) Let $M^n$ be an immersed submanifold of a real space form with constant sectional curvature $\kappa$. Then
\begin{equation}
\rho + \rho^\perp \leq ||H||^2 + \kappa,
\end{equation}
where $\rho$ is the normalized scalar curvature, an intrinsic invariant of the submanifold derived from the curvature tensor of the tangent bundle, $\rho^\perp$ is the normalized normal scalar curvature, an extrinsic invariant that captures the curvature of the normal bundle and $H$ 
is the normalized mean curvature vector field, which measures the average extrinsic curvature of the submanifold.

The Normal Scalar Curvature Conjecture was
 proved by Lu \cite{Lu} and Ge-Tang \cite{Ge}
 independently. The submersion version of the normal scalar
 curvature conjecture was studied by Ge \cite{Gee}.

The DDVV conjecture, as originally stated and proven, applies to submanifolds of real space forms. The problem of generalizing this inequality to more intricate ambient spaces and types of maps has become a significant area of contemporary research. The DDVV inequality has been successfully extended to submanifolds of warped products, as well as to special classes of submanifolds, such as slant, totally real, and CR-submanifolds, particularly in complex and quaternionic ambient spaces.

The extension to quaternionic space forms, which are the quaternionic analogs of real space forms, is particularly challenging. Unlike the real and complex numbers, the quaternions form a non-commutative division algebra. This non-commutativity is a fundamental property that distinguishes quaternionic geometry from its real and complex counterparts. 

 Recall that a Riemannian map was defined by Fischer \cite{Fish} to be a map whose derivative is a linear isometry between the domain tangent space modulo the kernel and its range. A Riemannian map is a smooth map between Riemannian manifolds that generalizes both isometric immersions and Riemannian submersions. Unlike an isometric immersion, which has an injective differential, a Riemannian map can have a non-trivial kernel (or vertical distribution) at each point.

Fischer presented in his foundational work \cite{Fish} that a Riemannian map “is a map that is as isometric as it can be, subject to the limitations imposed upon it as a differentiable mapping”. He demonstrated that these maps fulfill the generalized eikonal equation and thus have constant rank on each connected component,  establishing a connection between physical optics and geometric optics. Additionally, Fischer proposed a framework for developing a quantum model of nature using Riemannian maps, thereby linking the mathematical theories of Riemannian and harmonic maps with Lagrangian field theory, alongside the physical principles embodied in the Maxwell-Shr\"odinger equations.
 
The structure of the paper is outlined as follows. Section 2 revisits the fundamental concepts and formulas necessary for the discussion. In Section 3, we explore the DDVV inequality for Riemannian maps transitioning from a Riemannian manifold to a  quaternionic space form. Lastly, the paper concludes with the application of such inequalities.

\section{Preliminaries}
In this section, we review some basic concepts and results on geometric structures for Riemannian maps, focusing on the second fundamental form of a Riemannian map.
\subsection{Quaternionic K\"ahler Manifold}

Let $M$ be a $4m$-dimensional $C^{\infty}$-manifold and let $E$ be a rank 3
subbundle of End$(TM)$ such that for any point $p\in M$ with a neighborhood
$U$, there exists a local basis $\left\{ J_{1},J_{2},J_{3}\right\} $ of 
sections of $E$ on $U$ satisfying  
\[
J_{\alpha}^{2}=-id,\quad J_{\alpha }J_{\alpha +1}=-J_{\alpha +1}J_{\alpha
}=J_{\alpha +2}  
\]
where the indices $\alpha$ are taken from $\{1,2,3\}$
modulo $3$. Then $E$ is said to  be an almost quaternionic structure on $M$
and $(M,E)$ an almost  quaternionic manifold \cite{Alekseevsky-Marchiafava}.\\
Moreover, let $g_{1}$ be a Riemannian metric on $M$ defined by  
\[
g_{1}\left( J_{\alpha }X,J_{\alpha }Y\right) =g_{1}(X,Y), 
\]
for all vector fields $X,Y\in \Gamma (TM)$, where the indices $\alpha$ are taken from
$\{1,2,3\}$ modulo $3$. Then $(M,E,g_{1})$ is said to be an almost quaternionic Hermitian manifold \cite%
{Ianus-Mazzocco-Vilcu} and the basis $\left\{ J_{1},J_{2},J_{3}\right\} $ is
said to be a quaternionic Hermitian basis.

An almost quaternionic Hermitian manifold $(M,E,g)$ is said to be a 
quaternionic K\"{a}hler manifold \cite{Ianus-Mazzocco-Vilcu} if there exist 
locally defined $1$-forms $\omega _{1},\omega _{2},\omega _{3}$ such that 
for $\alpha \in \{1,2,3\}$ modulo $3$  
\begin{equation}
\nabla _{X}J_{\alpha }=\omega _{\alpha +2}(X)J_{\alpha+1}-\omega _{\alpha
+1}(X)J_{\alpha +2},  \label{eq-qkm}
\end{equation}
for $X\in \Gamma (TM)$, where the indices are taken from $\{1,2,3\}$ modulo $
3$.\\


\noindent Let $\left( M,E,g_{1}\right) $ be a quaternionic K\"ahler manifold and let $X$ be a non-null vector on $M.$
Then the $4$-plane spanned by $\{X,J_{1}X,J_{2}X,J_{3}X\}$, denoted by $Q(X),$  is called a quaternionic
$4$-plane. Any $2$-plane in $Q(X)$  is called a quaternionic plane. The sectional curvature of
a quaternionic plane is called a quaternionic sectional curvature. A quaternionic K\"ahler
manifold is a quaternionic space form if its quaternionic sectional curvatures are equal
to a constant, say $c.$ A quaternionic K\"ahler manifold $ (M, E, g_{1})$  is a
quaternionic space form, denoted $M(c),$  if and only if its Riemannian curvature tensor $R$ on $M(c)$ is given by \rm{\cite{S Ishi}}
\vspace{-.1cm}
\begin{equation} \label{eq-2.1}
\begin{aligned}
    R(X, Y)Z = \frac{c}{4} \Bigg\{ & g_1(Y, Z)X - g_1(X, Z)Y \\
    & + \Bigg[ \sum_{\alpha = 1}^{3} \Big(
        g_1(J_{\alpha}Y, Z) J_{\alpha}X 
        - g_1(J_{\alpha}X, Z) J_{\alpha}Y 
        + 2g_1(J_{\alpha}Y,X) J_{\alpha}Z
    \Big) \Bigg] \Bigg\},
\end{aligned}
\end{equation}
for $X,Y,Z\in \Gamma (TM)$.

\subsection{Riemannian Maps}
Consider $F: (M,g_{1}) \rightarrow (N,g_{2})$ be a smooth map between Riemannian manifolds $M$ and $N$ of dimension $m$ and $n,$ respectively, such that $0 < rank F < min \{m,n\}$ and if $F_{*}:T_{p}M \rightarrow T_{F(p)}N$ denotes the differential map at $p \in M,$ and $F(p) \in N,$ then $T_{p}M$ and $ T_{F(p)}N$ split orthogonally with respect to $g_{1}(p)$ and $g_{2}(F(p)),$ respectively, as \rm{\cite{Fish}}
\[ T_{p}M= \operatorname ker F_{*p} \oplus(\operatorname ker F_{*p})^{\perp} = {\mathcal{V}}_{p} \oplus{\mathcal{H}}_{p},\]
where ${\mathcal{V}}_{p}=\operatorname ker F_{*p}$ and ${\mathcal{H}}_{p}=(\operatorname ker F_{*p})^{\perp}$ are vertical and horizontal parts of $T_{p}M$ respectively.
Since $0 < rank F < min \{m,n\},$ we have $(range F_{*p})^\perp \neq 0.$ Therefore $T_{F(p)}N$ can be decomposed as follows: 
$$ T_{F(p)}N= {range F_{*}}_{p} \oplus (range F_{*p})^\perp. $$
Then the map $F: (M,g_{1}) \rightarrow (N,g_{2})$ is called a Riemannian map at $ p \in M,$ if $g_{2}(F_{*}X,F_{*}Y)= g_{1}(X,Y)$ for all vector fields $X,Y \in \Gamma \left(({\operatorname ker F_{*p}})^{\perp} \right).$

In particular, if $ker F_*=0$, then a Riemannian map is just an isometric immersion, while if $(range F_*)^\perp=0$
 then a Riemannian map is nothing but a Riemannian submersion.

 The second fundamental form of the map $F$ is given by \rm{\cite{Sahin17}}
\[ (\nabla F_{*})(X,Y)={\nabla}_{X}^{N}F_{*}Y-F_{*}({\nabla}_{X}^{M}Y),\]
where $\nabla ^{M}$ is the Levi-Civita connection on $M$ and $\nabla ^{N}$ is the pullback connection of $\nabla ^{M}$ along $F,$ provided that $(\nabla F_{*})(X,Y)$ has no components in $range F_{*},$ if $X,Y \in \Gamma \left(({\operatorname ker F_{*}})^{\perp} \right).$ \\



\noindent Let $\nabla ^{F\perp}$ be a linear connection on $ ({range F_{*}})^{\perp}$ and $ {\nabla }^{F\perp}_{X}V$ be the orthogonal projection of ${\nabla}_{X}^{N}V $ onto $(range F_{*})^{\perp}, X \in \Gamma({\operatorname(ker F_{*})^{\perp}}), V \in \Gamma \left(({range F_{*}})^{\perp} \right),$ then \rm{\cite{Sahin17}} 
\[ {\nabla}_{F_{*}X}^{N}V=-S_{V}F_{*}X+ {\nabla}_{X}^{F \perp}V,\]
where $-S_{V}F_{*}X$ is the tangential component of ${\nabla}_{F_{*}X}^{N}V$.
We have,
\begin{equation}\label{eq-6a}
    g_{2}\left( S_{V}F_{*}X,F_{*}Y\right)=g_{2}\left( V,(\nabla F_{*})(X,Y)\right),
\end{equation}
for $X,Y \in \Gamma \left(({\operatorname ker F_{*}})^{\perp} \right)$ and $ V \in \Gamma \left(({range F_{*}})^{\perp} \right).$ Since $\left(\nabla F_*\right)$ is symmetric, it follows that $S_V$ is a symmetric linear transformation of range $F_*$, called the shape operator of a Riemannian map $F$. $S_{V}F_{*}X$ is bilinear on $V$ and $F_{*}X$. \\
Let $R^M, R^N$ and $R^{\perp}$ be the curvature tensors of $\nabla^{M}, \nabla^{N}$ and ${(range F_{*})}^{\perp},$ respectively. Then we have the Gauss equation and the Ricci equation for the Riemannian map $F$ given by \rm{\cite{Sahin17}}
\begin{equation} \label{eq-2.2}
    \begin{aligned}
        g_{2}\left( R^{N} (F_{*}X,F_{*}Y)F_{*}Z, F_{*}H\right)&= g_{1}\left( R^{M} (X,Y)Z,H \right) + g_{2}\left( (\nabla F_{*})(X,Z),(\nabla F_{*})(Y,H)\right)\\
        & -g_{2}\left( (\nabla F_{*})(Y,Z), (\nabla F_{*})(X,H) \right),
    \end{aligned}
\end{equation}
\begin{equation} \label{eq-2.3}
    \begin{aligned}
        g_{2}\left( R^{N} (F_{*}X,F_{*}Y)V_{1}, V_{2}\right)&= g_{2}\left( R^{F \perp} (F_{*}X,F_{*}Y)V_{1},V_{2}\right) + g_{2}\left( [S_{V_{2}},S_{V_{1}}]F_{*}X,F_{*}Y\right),
    \end{aligned}
\end{equation}
where $X,Y,Z,H \in \Gamma({\operatorname(ker F_{*})^{\perp}}), V_{1},V_{2} \in \Gamma({(range F_{*})^{\perp}})$ and $[S_{V_{2}},S_{V_{1}}]=S_{V_{2}}S_{V_{1}}-S_{V_{1}}S_{V_{2}}.$
Let $\{e_{1},\ldots,e_{r}\} $ and $\{ v_{r+1},\ldots,v_{n}\}$ be the orthonormal basis of ${\operatorname(ker F_{*})^{\perp}}_{p} $ and ${(range F_{*p})^{\perp}}$, $p \in M,$ respectively. Then the Ricci curvature $ {\mathcal{R}}ic^{\operatorname(ker F_{*})^{\perp}}$, scalar curvature ${\tau^{\operatorname(ker F_{*})^{\perp}}}$ and normalized scalar curvature $ {\varrho^{\operatorname(ker F_{*})^{\perp}}}$ on $({\operatorname(ker F_{*})^{\perp}})_{p}$ are given by
\[ {\mathcal{R}}ic^{\operatorname(ker F_{*})^{\perp}}(X)=\sum_{i=1}^{r}g_{1}(R^{M}(e_{i},X)X,e_{i}), \quad \forall X \in {\operatorname(ker F_{*})^{\perp}},\]
\[\tau^{\operatorname(ker F_{*})^{\perp}}= \sum_{1\leq i<j \leq r}g_{1}(R^{M}(e_{i},e_{j})e_{j},e_{i}),\]
\[\varrho^{\operatorname(ker F_{*})^{\perp}}=\frac{2 \tau^{\operatorname(ker F_{*})^{\perp}}}{r(r-1)},\]
whereas, scalar curvature of sub-bundle $ (range F_{*})^{\perp}$ called normal scalar curvature, is defined as 
\[ \tau^{{(range F_{*})}^{\perp}}= \sqrt{\sum_{1 \leq i<j \leq r} \sum_{1 \leq \beta < \gamma \leq n}g_{2}(R^{\perp}(e_{i},e_{j})v_{\beta},v_{\gamma})^{2}},\] 
and normalized normal scalar curvature $ \varrho^{{(range F_{*})}^{\perp}}$ is given by 
\begin{equation}\label{eq-5a}
\varrho^{{(range F_{*})}^{\perp}}= \frac{\tau^{{(range F_{*})}^{\perp}}}{r(r-1)}.
\end{equation}
Also, for the Riemannian map $F,$ we can write \rm{\cite{Lee}}
\begin{equation}\label{eq-5}
 \zeta_{ij}^{\beta}= g_{2}\left( (\nabla F_{*})(e_{i},e_{j}),v_{\beta}\right),\quad  i,j=1 \ldots,r, \quad \beta=r+1, \ldots,n  
\end{equation}
\begin{equation}\label{eq-6}
    \lVert \zeta \rVert^{2}=\sum_{i,j=1}^{r}g_{2}\left( (\nabla F_{*})(e_{i},e_{j}),(\nabla F_{*})(e_{i},e_{j})\right),
\end{equation}
\begin{equation}\label{eq-7}
    trace \zeta =\sum_{i=1}^{r}(\nabla F_{*})(e_{i},e_{j}),
\end{equation}
\begin{equation}\label{eq-8}
    \lVert trace \zeta \rVert^{2}=g_{2}(trace \zeta, trace \zeta).
\end{equation}
Therefore, we have 
\begin{equation}\label{eq-62}
   {\mathcal C}^{\mathcal H}= \frac{1}{r}\lVert \zeta \rVert^{2},
\end{equation}
which is called the Casorati curvature of the horizontal space ${\mathcal H}^p$.

\begin{lem} \label{lem-1}
{\rm{\cite{Lee,Gulbahar}}} For the above terms, we have
\begin{equation}\label{eq-3.6} 
     \begin{aligned}
       r{\mathcal C}^{\mathcal H} = \lVert \zeta \rVert^{2}&= \frac{1}{2}\lVert trace \zeta \rVert^{2} + \frac{1}{2} \sum_{\beta=1}^{n}({\zeta}_{11}^{\beta}-{\zeta}_{22}^{\beta}- \ldots - {\zeta}_{rr}^{\beta})^{2}+ 2 \sum_{\beta=r+1}^{n} \sum_{i=2}^{r} ({\zeta}_{1i}^{\beta})^{2} \\
         & -2 \sum_{\beta=r+1}^{n} \sum_{2 \leq i < j \leq r} \Big\{ {\zeta}_{ii}^{\beta}{\zeta}_{jj}^{\beta}-\left({\zeta}_{ij}^{\beta}\right)^{2} \Big\}.
\end{aligned}
 \end{equation}
 \end{lem}
 \begin{lem} \label{lem-2}
     {\rm{\cite{Lu}}} Let $({\zeta}_{ij}^{\beta})$,  $i,j=1,\ldots,r$ and $\beta=1,\ldots,n$ be the entries of
 (the traceless part of) the second fundamental form under the orthonormal frames of both the
 tangent bundle and the normal bundle. Then 
\begin{equation}
    \begin{aligned}
       \sum_{\beta=1}^{n}\sum_{1 \leq i < j \leq r} ({\zeta}_{ii}^{\beta}-{\zeta}_{jj}^{\beta})^{2} + & 2r \sum_{\beta=1}^{n}\sum_{1 \leq i < j \leq r} ({\zeta}_{ij}^{\beta})^{2}\\ & \geq 2r \left( \sum_{1 \leq \beta < \gamma \leq n} \sum_{1 \leq i<j \leq r} \left( \sum_{k=1}^{r} ({\zeta}_{jk}^{\beta}{\zeta}_{ik}^{\gamma}-{\zeta}_{ik}^{\beta}{\zeta}_{jk}^{\gamma})\right)^{2} \right)^{\frac{1}{2}}.  
    \end{aligned}
\end{equation}
 \end{lem}
\section{Riemannian Maps to Quaternionic Space Forms}

Let $\{e_{1},\ldots, e_{r}\}$ and $ \{v_{r+1},\ldots, v_{n}\}$ be the orthonormal basis of ${{(\operatorname ker F_{*})}^{\perp}}_{p}$ and $ {{(range F_{*})}^{\perp}}_{F_{*}p},$ $ p \in M,$ respectively, then $ \{F_{*}e_{1}, \ldots, F_{*}e_{r}\}$ is the orthonormal basis of $ {{(range F_{*})}}_{F_{*}p}.$
\begin{thm}
Let $F:(M,g_{1}) \rightarrow (M(c),g_{2})$ be a Riemannian map from a Riemannian manifold $M$ to quaternionic space form $M(c)$ and rank $F= r< n$ then we have 
\begin{equation} \label{eq-1}
4{\mathcal{R}}ic^{\operatorname(ker F_{*})^{\perp}}(X) \leq c (r-1)+ \lVert trace \zeta \rVert ^{2} + 3c \sum_{\alpha = 1}^{3} \sum_{i = 1}^{r} g_{2}^{2}\left( J_{\alpha}F_{*}X,F_{*}e_{i}\right),
\end{equation}
where X is a unit horizontal vector field on $M$. The equality case of {\rm(\ref{eq-1})} arises for $X \in \Gamma ((\operatorname{ker F_{*}})^{\perp})$ if and only if 
    \begin{equation}
    \label{eq-3.2}(\nabla F_{*})(X,Y)=0 \quad \forall \ Y \in \Gamma ((\operatorname{ker F_{*}})^{\perp}) \text{ orthogonal to} \ X, 
    \end{equation}
    \begin{equation}\label{eq-3.3}
   {\text and}\quad  (\nabla F_{*})(X,X)= \frac{1}{2} \operatorname{trace} \zeta. 
    \end{equation}
    \end{thm}
{\bf Proof:} Let $M(c)$ be a quaternionic space form with constant sectional curvature $c,$ then from \rm(\ref{eq-2.2}) and \rm(\ref{eq-2.3}), we have 
\begin{equation}\label{eq-3.4}
    \begin{aligned}
        g_{1}\left( R^{M}(X,Y)Z,H\right) &=\frac{c}{4} \Bigg\{ g_{1}(Y, Z)g_{1}(X,H) - g_{1}(X, Z)g_{1}(Y,H) 
     + \Bigg[ \sum_{\alpha = 1}^{3} \Big(
        g_{2}(J_{\alpha}F_{*}Y,F_{*}Z) g_{2}(J_{\alpha}F_{*}X,F_{*}H)\\ 
        &- g_{2}(J_{\alpha}F_{*}X,F_{*} Z) g_{2}(J_{\alpha}F_{*}Y,F_{*}H)
        + 2g_{2}(J_{\alpha}F_{*}Y,F_{*} X) g_{2}(J_{\alpha}F_{*}Z,F_{*}H)
    \Big) \Bigg] \Bigg\} \\ &-g_{2}\left( (\nabla F_{*})(X,Z), (\nabla F_{*})(Y,H)\right) + g_{2}\left( (\nabla F_{*})(Y,Z), (\nabla F_{*})(X,H)\right),
    \end{aligned}
\end{equation}

\noindent where $X,Y,Z,H \in \Gamma ((\operatorname{ker F_{*}})^{\perp}).$ Let$\{e_{1}, \ldots, e_{r}\}$ with $e_{1}=X$ and using \rm(\ref{eq-5})-\rm(\ref{eq-8}) and \rm(\ref{eq-3.4}), we get
 \begin{equation} \label{eq-3.5}
  \frac{c}{4}r(r-1) + \frac{c}{4} \sum_{\alpha=1}^{3} \sum_{i,j=1}^{r} \Big[ 3g_{2}^{2} (J_{\alpha}F_{*}e_{i},F_{*}e_{j}) \Big] = \lVert \zeta \rVert ^{2}- \lVert trace \zeta \rVert ^{2} + 2 \tau ^{{(\operatorname ker F_{*})}^\perp}.
 \end{equation}
 using lemma (\ref{lem-1}), and 
\begin{equation} \label{eq-3.7}
    \begin{aligned}
      \tau ^{{(\operatorname ker F_{*})}^\perp} (p)&= \frac{c}{8}(r-1)(r-2)+ \frac{c}{4}\sum_{\alpha=1}^{3} \Big[ 3 \sum_{2 \leq i < j \leq r}g_{2}^{2} (J_{\alpha}F_{*}e_{i},F_{*}e_{j})\Big] \\ &- \sum_{\beta=r+1}^{n} \sum_{2 \leq i < j \leq r} \Big\{\left({\zeta}_{ij}^{\beta}\right)^{2} - {\zeta}_{ii}^{\beta}{\zeta}_{jj}^{\beta} \Big\} + {\mathcal{R}}ic^{{{(\operatorname ker F_{*})}^\perp}}(X).
    \end{aligned}
\end{equation}
Substituting the value of \rm(\ref{eq-3.6}) and \rm(\ref{eq-3.7}) in \rm(\ref{eq-3.5}), we get 
\begin{equation}
    \begin{aligned}
{\mathcal{R}}ic^{{{(\operatorname ker F_{*})}^\perp}}(X)&=\frac{c}{4}(r-1)+ \frac{1}{4}\lVert trace \zeta \rVert ^{2} + \frac{3c}{4} \sum_{\alpha = 1}^{3} \sum_{i = 1}^{r} g_{2}^{2}\left( J_{\alpha}F_{*}X,F_{*}e_{i}\right)\\ &-\frac{1}{4} \sum_{\beta=r+1}^{n}({\zeta}_{11}^{\beta}-{\zeta}_{22}^{\beta}- \ldots - {\zeta}_{rr}^{\beta})^{2}- \sum_{\beta=r+1}^{n} \sum_{i=2}^{r} ({\zeta}_{1i}^{\beta})^{2}\\ & \leq \frac{c}{4}(r-1)+ \frac{1}{4}\lVert trace \zeta \rVert ^{2} + \frac{3c}{4} \sum_{\alpha = 1}^{3} \sum_{i = 1}^{r} g_{2}^{2}\left( J_{\alpha}F_{*}X,F_{*}e_{i}\right).
    \end{aligned}
\end{equation}
Hence, we have the inequality \rm(\ref{eq-1}).
The equality case of \rm(\ref{eq-1}) holds if and only if 
\[ {\zeta}_{11}^{\beta}={\zeta}_{22}^{\beta}+ \ldots +{\zeta}_{rr}^{\beta}\]
and \[{\zeta}_{1i}^{\beta}=0, \quad i=2, \ldots,r;\quad\beta=r+1, \ldots, n;\] hence we have \rm(\ref{eq-3.2}).

\begin{thm}
Let $F: (M,g_{1}) \rightarrow (M(c),g_{2})$ be a Riemannian map from a Riemannian manifold $M$ to quaternionic space form $M(c)$ with rank$F=r<n,$ then normalized normal scalar curvature $\varrho ^{{(range F_{*})}^{\perp}}$ and normalized scalar curvature  $\varrho ^{{(\operatorname ker F_{*})}^{\perp}}$ satisfy
\begin{equation} \label{eq-3.9}
 \varrho ^{{(range F_{*})}^{\perp}} + \varrho ^{{(\operatorname ker F_{*})}^{\perp}} \leq \frac{1}{r^{2}} \lVert trace \zeta \rVert^{2}+ \frac{c}{4} + \frac{3c}{4r(r-1)}\sum_{\alpha=1}^{3}\sum_{i,j=1}^{r}g_{2}^{2}\left( J_{\alpha}F_{*}e_{i},F_{*}e_{j}\right). 
\end{equation}
\end{thm}
{\bf Proof:} Let $\{e_{1},\ldots, e_{r}\}$ and $ \{v_{r+1},\ldots, v_{n}\}$ be the orthonormal basis of ${{(\operatorname ker F_{*})}^{\perp}}_{p}$ and $ {{(range F_{*})}^{\perp}}_{F_{*}p},$ $ p \in M,$ respectively, then $ \{F_{*}e_{1}, \ldots, F_{*}e_{r}\}$ is the orthonormal basis of $ {{(range F_{*})}}_{F_{*}p}.$ Then from \rm(\ref{eq-2.1}),\rm(\ref{eq-2.2}) and \rm(\ref{eq-2.3}), we have
\begin{equation} \label{eq-3.10}
    \begin{aligned}
\sum_{\beta,\gamma=r+1}^{n}\sum_{i,j=1}^{r} g_{2}\left( R^{F \perp} (F_{*}e_{i},F_{*}e_{j})v_{\beta},v_{\gamma}\right)&=  g_{2}\left( R^{N} (F_{*}e_{i},F_{*}e_{j})v_{\beta}, v_{\gamma}\right)- g_{2}\left( [S_{v_{\gamma}},S_{v_{\beta}}]F_{*}e_{i},F_{*}e_{j}\right),
    \end{aligned}
\end{equation}
using \rm(\ref{eq-3.10}) in \rm(\ref{eq-5a}), we get
\begin{equation}\label{eq-5b}
\varrho^{{(range F_{*})}^{\perp}}= \frac{2}{r(r-1)}\sqrt{\sum_{1 \leq i<j \leq r} \sum_{1 \leq \beta < \gamma \leq n} g_{2}^{2}\left( [S_{v_{\gamma}},S_{v_{\beta}}]F_{*}e_{i},F_{*}e_{j}\right)},
\end{equation}
using \rm(\ref{eq-6a}) in \rm(\ref{eq-5b}) and solving further, we get 
\begin{equation}\label{eq-5c}
\varrho^{{(range F_{*})}^{\perp}}= \frac{2}{r(r-1)}\sqrt{\sum_{1 \leq i<j \leq r} \sum_{1 \leq \beta < \gamma \leq n} \left( \sum_{k=1}^{r} ({\zeta}_{jk}^{\beta}{\zeta}_{ik}^{\gamma}-{\zeta}_{ik}^{\beta}{\zeta}_{jk}^{\gamma})\right)^{2}},
\end{equation}
using lemma (\ref{lem-2}), in \rm(\ref{eq-5c}), we get

\begin{equation}\label{eq-5d}
    \begin{aligned}
        r^{2}(r-1)\varrho^{{(range F_{*})}^{\perp}} \leq \sum_{\beta=r+1}^{n}\sum_{1 \leq i < j \leq r} ({\zeta}_{ii}^{\beta}-{\zeta}_{jj}^{\beta})^{2} + 2r \sum_{\beta=r+1}^{n}\sum_{1 \leq i < j \leq r}  ({\zeta}_{ij}^{\beta})^{2}.
    \end{aligned}
\end{equation}
(Remark- ${\zeta}_{ii}^{\beta}=0$ for $\beta= 1,\ldots, r $, from equation {\rm(\ref{eq-5})})\\
Also, we can compute
\begin{equation}\label{eq-5e}
    (r-1)\lVert trace \zeta\rVert ^{2}=\sum_{\beta=r+1}^{n}\sum_{1 \leq i < j \leq r}\left(({\zeta}_{ii}^{\beta}-{\zeta}_{jj}^{\beta})^{2} + 2r {\zeta}_{ii}^{\beta}{\zeta}_{jj}^{\beta} \right), 
\end{equation}
and 
\begin{equation}\label{eq-5f}
    \begin{aligned}
     \tau ^{{(\operatorname ker F_{*})}^{\perp}}= \frac{r-1}{4}\left( \frac{c}{2}r\right) + \frac{3c}{8}\sum_{\alpha=1}^{3}\sum_{i,j=1}^{r}g_{2}^{2}(J_{\alpha}F_{*}e_{i},F_{*}e_{j})+ \sum_{\beta=r+1}^{n}\sum_{1 \leq i < j \leq r}\left( {\zeta}_{ii}^{\beta}{\zeta}_{jj}^{\beta}-({\zeta}_{ij}^{\beta})^{2}\right).  
    \end{aligned}
\end{equation}
From \rm(\ref{eq-5d}), \rm(\ref{eq-5e}) and \rm(\ref{eq-5f}), we get the required result.

\medskip

\noindent{\bf Remark:} The establishment of a DDVV-type inequality for Riemannian maps from a Riemannian manifold into a quaternionic space form represents a significant advancement in differential geometry. The core application of the DDVV inequality for Riemannian maps is to establish a relationship between the curvature of the source manifold and the extrinsic curvature of the map. They serve as a fundamental way to measure how a map "bends" or "curves" as it projects one manifold onto another. When the inequality becomes an equality, it characterizes a specific class of maps with special geometric properties. These maps are considered ``ideal" or ``optimal" in some geometric sense and can lead to rigidity results, which state that a map satisfying the equality must be of a very specific, rigid type. This work provides a unifying framework that bridges several distinct and complex areas of research: the DDVV conjecture, the non-commutative geometry of quaternions, and the generalized mapping theory of Riemannian maps.

The results of this paper can be applied to a variety of related fields. In mathematical physics, quaternionic K\"ahler manifolds are relevant in the study of supersymmetry, string theory, and general relativity. Riemannian maps and submersions are used to model phenomena in Kaluza-Klein theory, where higher-dimensional manifolds are projected onto lower-dimensional ones. The DDVV-type inequality can provide new analytical tools for these applications, offering constraints and relationships between physical quantities that are geometrically motivated.

\medskip 
\noindent \textbf {Declarations:}\\
\textbf { Data Availability:} No new data were created or analysed in this study.\\
\textbf {Conflicts of interest:} The authors declare that they have no known financial conflicts of
interest or personal relationships that could have influenced the work presented in this paper.\\
\textbf {Funding:} The corresponding author is supported and funded by the National Board of Higher Mathematics (NBHM) project no. 02011/21/2023NBHM(R.P.)/R\&DII/14960, INDIA.\\
\textbf {Ethics approval:} The authors hereby affirm that the contents of this manuscript are original.
Furthermore, it has neither been published elsewhere in any language fully or partly, nor is it under
review for publication anywhere.


\begin{thebibliography}{99}

\bibitem{Alekseevsky-Marchiafava} D.V. Alekseevsky and S. Marchiafava: 
Almost complex submanifolds of quaternionic manifolds, Steps in Diff. Geom., Institutes of Mathematics and Informatics, 
University of Debrecen, (2001), 23-38.



 


 


 \bibitem{DDVV} P.J. De Smet, F. Dillen, Leopold C.A. Verstraelen and L. Vrancken: A pointwise inequality in submanifold theory, Arch. Math., 35 (1999), 115-128.



\bibitem{Fish} A.E. Fischer: Riemannian Maps Between Riemannian Manifolds. Mathematical Aspects of Classical Field Theory (Seattle, WA, 1991). Contemporary Mathematics, American Mathematical Society, Providence, 132 (1992), 331-366.

\bibitem{Gee} J. Ge: DDVV-type inequality for skew-symmetric matrices and Simons-type inequality for Riemannian submersions,
Adv. Math., 251 (2014), 62-86.

\bibitem{Ge} J. Ge and Z. Tang: A proof of the DDVV conjecture and its equality case, Pacific Jour. Math., 237 (2008), no. 1, 87-95.

\bibitem{Guad} I.V. Guadalupe and L. Rodriguez: Normal curvature of surfaces in space forms, Pacific Jour. Math., 106 (1983), 95-103.


\bibitem{Gulbahar} M. G\"ulbahar, \c S.E. Meri\c c and E. Kili\c c: Sharp inequalities involving the Ricci curvature for Riemannian submersions, Kragujevac Jour. Math., 41 (2017), no. 2, 279-293.

\bibitem{Ianus-Mazzocco-Vilcu} S. Ianus, R. Mazzocco and G.E. Vilcu: Riemannian submersions from quaternionic manifolds, Acta. Appl. Math., 104 (2008), 83-89. 

\bibitem{S Ishi} S. Ishihara: Quaternion K\"ahlerian manifolds, Jour. Diff. Geom., 9 (1974), 483-500.

\bibitem{Lu} Z. Lu: Normal scalar curvature conjecture and its applications, Jour. Funct. Anal., 261 (2011), no. 5, 1284-1308.

\bibitem{Lee} J.W. Lee, C.W. Lee, B. \c Sahin and G.E. Vilcu: Chen-Ricci inequalities for Riemannian maps and their applications, Contemp. Math., 777 (2022), 137-152.








\bibitem{Sahin17}  B. \c Sahin: Riemannian Submersions, Riemannian maps in Hermitian Geometry and their Applications, Cambridge, Elsevier Academic Press, (2017).

 

\bibitem{Wintgen} P. Wintgen: Sur l’in\'egalit\'e de Chen-Willmore, C. R. Acad. Sci., Paris, S\'er. A, 288  (1979), 993-995.
\end{thebibliography}
\end{document}